\documentclass[aop,MSNbibl,seceqn,citesort,dvips]{arximspdf}
\usepackage{graphicx}

\doi{10.1214/11-AOP664}
\volume{39}
\issue{5}
\pubyear{2011}
\firstpage{2018}
\lastpage{2041}

\makeatletter

\def\bptnote#1{}

\def\E{{\mathbb E}}
\def\N{{\mathbb N}}
\def\P{{\mathbb P}}
\def\R{{\mathbb R}}
\def\Z{{\mathbb Z}}
\def\B{{\mathbb B}}

\def\cE{{\mathcal E}}
\def\cB{{\mathcal B}}

\def\eps{\varepsilon}

\def\bk{\setminus}
\def\ind{{1}}

\newtheorem{theorem}{Theorem}[section]
 \newtheorem{lemma}[theorem]{Lemma}
 \newtheorem{proposition}[theorem]{Proposition}
 \newproclaim{definition}[theorem]{Definition}
 \newtheorem{conjecture}[theorem]{Conjecture}

\makeatother

\begin{document}
\begin{frontmatter}

\title{A crossover for the bad configurations of random walk in random scenery}
\runtitle{A crossover for the bad configurations of random walk}

\begin{aug}
\author[A]{\fnms{S\'{e}bastien} \snm{Blach\`{e}re}\thanksref{a1,a2}\ead[label=e1]{s.blachere@gmail.com}},
\author[B]{\fnms{Frank} \snm{den Hollander}\thanksref{a2}\ead[label=e2]{denholla@math.leidenuniv.nl}}\\
\and
\author[C]{\fnms{Jeffrey E.} \snm{Steif}\corref{}\thanksref{a3}\ead[label=e3]{steif@chalmers.se}}
\runauthor{S. Blach\`{e}re, F. den Hollander and J. E. Steif}
\affiliation{LATP, University Aix-Marseille 1 and EURANDOM, Leiden
University and EURANDOM, and Chalmers University of Technology
 and G\"{o}teborg University}
\dedicated{This paper is dedicated to the memory of Oded Schramm}
\address[A]{S. Blach\`{e}re\\
LATP\\
University Aix-Marseille 1\\ 39 rue Joliot Curie, 13453 Marseille
Cedex\\ France\\
and\\
EURANDOM\\
P.O. Box 513, 5600 MB Eindhoven\\
The Netherlands\\
\printead{e1}} 
\address[B]{F. den Hollander\\
Mathematical Institute\\ Leiden University\\ P.O. Box 9512,
2300 RA Leiden\\ The Netherlands\\
and\\
EURANDOM\\
P.O. Box 513, 5600 MB Eindhoven\\
The Netherlands\\
\printead{e2}}
\address[C]{J. E. Steif\\Mathematical Sciences\\
Chalmers University of Technology\\
and\\
G\"{o}teborg University\\
SE-41296 Gothenburg\\ Sweden\\
\printead{e3}}
\end{aug}
\thankstext{a1}{Supported in part by EURANDOM in Eindhoven.}
\thankstext{a2}{Supported in part by DFG and NWO through the
Dutch-German Bilateral Research Group on ``Mathematics of Random
Spatial Models
from Physics and Biology'' (2004--2009).}
\thankstext{a3}{Supported in part by the Swedish Research Council and
by the G\"{o}ran Gustafsson
Foundation for Research in the Natural Sciences and Medicine.}
\received{\smonth{1} \syear{2010}}
\revised{\smonth{3} \syear{2011}}

%
\begin{abstract}
In this paper, we consider a random walk and a random color scenery on~$\Z$.
The increments of the walk and the colors of the scenery are
assumed to be i.i.d. and to be independent of each other. We are
interested in the random process of
colors seen by the walk in the course of time. Bad configurations for
this random
process are the discontinuity points of the conditional probability distribution
for the color seen at time zero given the colors seen at all later times.

We focus on the case where the random walk has increments $0$, $+1$ or
$-1$ with
probability $\eps$, $(1-\eps)p$ and $(1-\eps)(1-p)$, respectively,
with $p \in
[\frac12,1]$ and $\eps\in[0,1)$, and where the scenery assigns the
color black or
white to the sites of $\Z$ with probability $\frac12$ each. We show
that, remarkably,
the set of bad configurations exhibits a crossover: for $\eps=0$ and
$p \in(\frac12,
\frac45)$ all configurations are bad, while for $(p,\eps)$ in an
open neighborhood
of $(1,0)$ all configurations are good. In addition, we show that for
$\eps=0$ and
$p=\frac12$ both bad and good configurations exist. We conjecture that
for all $\eps
\in[0,1)$ the crossover value is unique and equals $\frac45$.
Finally, we suggest
an approach to handle the seemingly more difficult case where $\eps>0$
and $p\in
[\frac12,\frac45)$, which will be pursued in future work.
\end{abstract}

%

\end{frontmatter}

\section{Introduction}
\label{S1}


\subsection{Random walk in random scenery}
\label{S1.1}

We begin by defining the random process that will be the object of our study.
Let $X=(X_n)_{n \in\N}$ be i.i.d. random variables taking the values
$0$, $+1$
and $-1$ with probability $\eps$, $p(1-\eps)$ and $(1-p)(1-\eps)$,
respectively,
with $\eps\in[0,1)$ and $p\in[\frac12,1]$. Let $S=(S_n)_{n \in\N
_0}$ with
$\N_0:=\N\cup\{0\}$ be the corresponding \textit{random walk} on $\Z
$, defined
by
\[
S_0 := 0 \quad\mbox{and} \quad S_n := X_1 + \cdots + X_n, \qquad  n\in\N,
\]
that is, $X_n$ is the step at time $n$ and $S_n$ is the position at
time $n$. Let
$C=(C_z)_{z \in\Z}$ be i.i.d. random variables taking the values $B$ (black)
and $W$ (white) with probability $\frac12$ each. We will refer to $C$
as the
\textit{random coloring} of $\Z$, that is, $C_z$ is the color of site
$z$. The pair
$(S,C)$ is referred to as the \textit{random walk in random scenery} associated
with $X$ and $C$.

Let
\[
Y := (Y_n)_{n \in\N_0} \qquad\mbox{where }  Y_n := C_{S_n}
\]
be the sequence of colors observed along the walk. We will refer to $Y$
as the
\textit{random color record}. This random process, which takes values in
the set
$\Omega_0=\{B,W\}^{\N_0}$ and has full support on $\Omega_0$, will
be our main
object of study. Because the walk may return to sites it has visited
before and
see the same color, $Y$ has intricate dependencies. An overview of the ergodic
properties of $Y$ is given in \cite{HS05}.

We will use the symbol $\P$ to denote the joint probability law of $X$ and~$C$.
The question that we will address in this paper is whether or not there
exists a version $V(B\mid\eta)$ of the conditional probability
\[
\P (Y_0=B \mid Y=\eta\mbox{ on } \N ),
\qquad\eta\in\Omega_0,
\]
such that the map $\eta\mapsto V(B\mid\eta)$ is everywhere
continuous on
$\Omega_0$. It will turn out that the answer depends on the choice of
$p$ and
$\eps$.

In \cite{HSW04}, we considered the pair $(X,Y)$ and identified the
structure of
the set of points of discontinuity for the analogue of the conditional
probability
in the last display. However, $(X,Y)$ is much easier to analyze than
$Y$, because
knowledge of $X$ and $Y$ fixes the coloring on the support of $X$. Consequently,
the structure of the set of points of discontinuity for $(X,Y)$ is very
different
from that for $Y$. The same continuity question arises for the
two-sided version of
$Y$ where time is indexed by~$\Z$, that is, the random walk is
extended to negative
times by putting $S_0=0$ and $S_n-S_{n-1}=X_n$, $n\in\Z$, with $X_n$
the step at
time $n\in\Z$. In the present paper, we will restrict ourselves to
the one-sided version.

The continuity question has been addressed in the literature for a
variety of
random processes. Typical examples include Gibbs random fields that are
subjected
to some transformation, such as projection onto a lower-dimensional
subspace or evolution
under a random dynamics. It turns out that even simple transformations
can create\vadjust{\eject}
discontinuities and thereby destroy the Gibbs property. For a recent overview,
see \cite{ELR04}. Our main result, described in Section \ref{S1.4}
below, is a
contribution to this area.


\subsection{Bad configurations and discontinuity points}
\label{S1.2}

In this section, we view the conditional probability distribution of
$Y_0$ given
$(Y_n)_{n\in\N}$ as a map from $\Omega= \{B,W\}^{\N}$ to the set of
probability
measures on $\{B,W\}$ (as opposed to a map from $\Omega_0$ to this
set).
Our question about continuity of conditional probabilities will be
formulated in terms of so-called \textit{bad configurations}.

\begin{definition}\label{badconf}
Let $\P$ denote any probability measure on $\Omega_0$ with full
support. A configuration $\eta\in\Omega$ is said to be a bad
configuration if there is a $\delta>0$ such
that for all $m \in\N$ there are $n\in\N$ and $\zeta\in\Omega$,
with $n > m$ and $\zeta=\eta$ on $(0,m)\cap\N$, such that
\[
 \bigl|\P \bigl(Y_0=B \mid Y=\eta\mbox{ on } (0,n)\cap\N \bigr)
- \P \bigl(Y_0=B \mid Y=\zeta\mbox{ on } (0,n)\cap\N \bigr) \bigr|
\geq\delta.
\]
\end{definition}

In words, a configuration $\eta$ is bad when, no matter how large we
take~$m$, by tampering with $\eta$ inside $[m,n)\cap\N$ for some $n>m$ while
keeping it fixed
inside $(0,m)\cap\N$, we can affect the conditional probability
distribution of~$Y_0$ in a nontrivial way.
Typically, $\delta$ depends on $\eta$, while $n$ depends on $m$.
A~configuration that is not bad is called a \textit{good configuration}.

The bad configurations are the discontinuity points of the conditional
probability
distribution of $Y_0$, as made precise by the following proposition (see
\cite{MRV99}, Proposition 6, and \cite{HSW04}, Theorem 1.2).

\begin{proposition}\label{bad=disc}
Let $\B$ denote the set of bad configurations for $Y_0$.
\begin{longlist}[(ii)]
\item[(i)] For any version $V(B\mid\eta)$ of the conditional probability
$\P(Y_0=B\mid Y=\eta\mbox{ on } \N)$, the set $\B$ is contained in
the set of discontinuity points for the map $\eta\mapsto V(B\mid\eta
)$.
\item[(ii)] There is a version $V(B\mid\eta)$ of the conditional probability
$\P(Y_0=B\mid Y=\eta\mbox{ on } \N)$ such that $\B$ is equal to the
set of discontinuity points for the map $\eta\mapsto V(B\mid\eta)$.
\end{longlist}
\end{proposition}


\subsection{An educated guess}
\label{S1.3}

For the random color record, a \textit{naive guess} is that all
configurations are bad
when $p=\frac12$ because the random walk is recurrent, while all
configurations are
good when $p \in(\frac12,1]$ because the random walk is transient.
Indeed, in the
recurrent case we obtain new information about $Y_0$ at infinitely many times,
corresponding to the return times of the random walk to the origin,
while in the
transient case no such information is obtained after a finite time.
However, we will
see that this naive guess is wrong. Before we state our main result,
let us make an
\textit{educated guess}:
\begin{itemize}
\item
(EG1) $\forall  p \in[\frac12,\frac45]   \ \forall  \eps\in[0,1)
\dvtx   \B=\Omega$.
\item
(EG2) $\forall  p \in(\frac45,1]   \ \forall  \eps\in[0,1)
\dvtx   \B=\varnothing$.
\end{itemize}

The explanation behind this is as follows.

\textit{Fully biased.} Suppose that $p=1$. Then
\[
\P (Y_0=Y_1 \mid Y=\eta\mbox{ on } \N )= \eps+(1-\eps
)\tfrac12,
\]
where we use that, for any $p$ and $\varepsilon$, $S_1$ and $(Y_n)_{n\in
\N}$ are
independent. Hence, the color seen at time $0$ only depends on the
color seen
at time $1$, so that $\B=\varnothing$. (Note that if $\eps=0$, then
$Y$ is i.i.d.)

\textit{Monotonicity.} For fixed $\eps$, we expect monotonicity in $p$:
if a configuration is bad for some $p \in(\frac12,1)$, then it should
be bad for all $p' \in[\frac12,p)$ also. Intuitively, the random
walk with
parameters $(p',\eps)$ is exponentially more likely to return to $0$ after
time $m$ than the random walk with parameters~$(p,\eps)$, and therefore
we expect that it is easier to affect the color at $0$ for $(p',\eps)$
than for $(p,\eps)$.

\textit{Critical value.} For a configuration to be good, we expect that
the random
walk must have a strictly positive speed conditional on the color
record. Indeed,
only then do we expect that it is exponentially unlikely to influence
the color at
$0$ by changing the color record after time $m$. To compute the
threshold value
for $p$ above which the random walk has a strictly positive speed, let
us consider
the monochromatic configuration ``all black.'' The probability for the
random walk
with parameters $(p,\eps)$ to behave up to time~$n$ like a random walk
with parameters
$(q,\delta)$, with $q \in[\frac12,1]$ and $\delta\in[0,1)$, is
\[
e^{-nH((q,\delta) \mid(p,\eps))},
\]
where
\begin{eqnarray*}
H ((q,\delta) \mid(p,\eps) )
&:=& \delta\log \biggl(\frac{\delta}{\eps} \biggr)
+ (1-\delta)\log \biggl(\frac{1-\delta}{1-\eps} \biggr)\\
&&{} + (1-\delta) \biggl[q \log \biggl(\frac{q}{p} \biggr)
+ (1-q) \log \biggl(\frac{1-q}{1-p} \biggr) \biggr]
\end{eqnarray*}
is the relative entropy of the step distribution $(q,\delta)$ with
respect to the step
distribution $(p,\eps)$. The probability for the random coloring to be
black all the
way up to site $(1-\delta) (2q-1)n$ is
\[
\bigl(\tfrac12\bigr)^{(1-\delta)(2q-1)n}.
\]
The total probability is therefore
\[
e^{-nC(q,\delta)} \qquad\mbox{with }
C(q,\delta) := H ((q,\delta) \mid(p,\eps) )+(1-\delta
)(2q-1)\log2.
\]
The question is: For fixed $(p,\eps)$ and $n\to\infty$, does the
lowest cost occur for
$q=\frac12$ or for $q>\frac12$? Now, it is easily checked that $q
\mapsto C(q,\delta)$
is strictly convex and has a derivative at $q=\frac12$ that is
strictly positive if and
only if $p \in[\frac12,\frac45)$, irrespective of the value of
$\eps$ and $\delta$.
Hence, zero drift has the lowest cost when $p \in[\frac12,\frac
45]$, while strictly
positive drift has the lowest cost when $p \in(\frac45,1]$. This
explains (EG1) and (EG2).


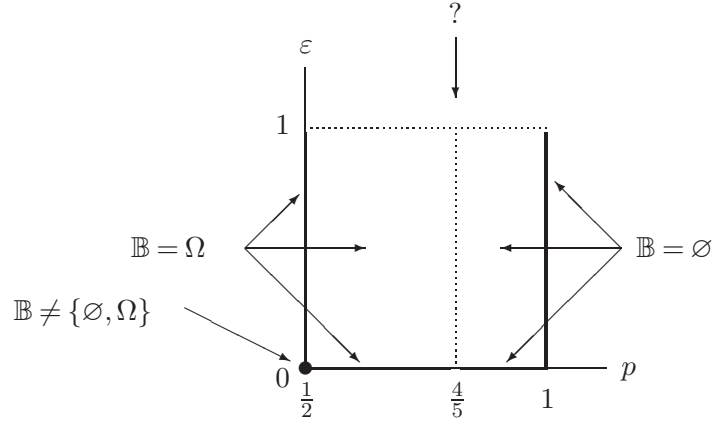
\begin{figure}
\begin{center}
\vspace*{1.8cm}
\setlength{\unitlength}{0.4cm}
\hspace*{2.4cm}
\begin{picture}(10,10)(0,-2)
\put(0,0){\line(10,0){10}}
\put(0,0){\line(0,10){10}}
\qbezier[40](5,0)(5,4)(5,8)
\qbezier[40](0,8)(4,8)(8,8)
\put(-2,4){\vector(1,0){4}}
\put(-2,4){\vector(1,-1){3.8}}
\put(-2,4){\vector(1,1){1.8}}
\put(10.5,4){\vector(-1,0){4}}
\put(10.5,4){\vector(-1,-1){3.8}}
\put(10.5,4){\vector(-1,1){2.2}}
\put(5,11){\vector(0,-1){2}}
\put(-4,2){\vector(2,-1){3.5}}
\put(11,3.7){$\B=\varnothing$}
\put(-5.8,3.7){$\B=\Omega$}
\put(-9.7,1.6){$\B\neq\{\varnothing,\Omega\}$}
\put(4.75,11.5){?}
\put(10.5,-0.2){$p$}
\put(-0.2,10.5){$\varepsilon$}
\put(-0.25,-1.3){$\frac12$}
\put(4.75,-1.3){$\frac45$}
\put(7.8,-1.3){$1$}
\put(-1,7.8){$1$}
\put(-1,-0.6){$0$}
\put(0,0){\circle*{0.4}}
{\thicklines
\qbezier(8,0)(8,4)(8,7.8)}
{\thicklines
\qbezier(0,0)(2,0)(4.8,0)}
{\thicklines
\qbezier(5.2,0)(7,0)(8,0)}
{\thicklines
\qbezier(0,0)(0,4)(0,7.8)}
\end{picture}
\end{center}
\caption{Conjectured behavior of the set $\B$ as a
function of $p$
and $\varepsilon$. Theorem \protect\ref{phtr} proves this behavior on the left
part of the
bottom horizontal line and in a neighborhood of the bottom right corner.}
\label{figbad}
\end{figure}


\subsection{Main theorem}
\label{S1.4}

We are now ready to state our main result and compare it with the
educated guess
made in Section \ref{S1.3} (see Figure~\ref{figbad}).

\begin{theorem}\label{phtr}
(\textup{i}) There exists a neighborhood of $(1,0)$ in the $(p,\varepsilon)$-plane for
which $\B=\varnothing$. This neighborhood can be taken to contain
the line segment $(p_*,1]
\times\{0\}$ with $p_* = 1/(1+5^5 12^{-6}) \approx0.997$.

(\textup{ii}) If $p \in(\frac12,\frac45)$ and $\eps=0$, then $\B=\Omega
$.

(\textup{iii}) If $p=\frac12$ and $\eps=0$, then $\B\notin\{\varnothing
,\Omega\}$.
\end{theorem}

Theorem \ref{phtr}(ii) and (iii) prove (EG1) for $p \in
[\frac12,\frac45)$ and $\eps=0$,
except for $p=\frac12$ and $\eps=0$, where (EG1) fails. We will see
that this
failure comes from parity restrictions. Theorem \ref{phtr}(i) proves
(EG2) in a
neighborhood of $(1,0)$ in the $(p,\varepsilon)$-plane. We already have
seen that
$\B=\varnothing$ when $p=1$ and $\eps\in[0,1)$. Note that Theorem
\ref{phtr}(ii) and (iii)
disprove monotonicity in $p$ for $\eps=0$. We believe this
monotonicity to fail
only at $p=\frac12$ and $\eps=0$.

To appreciate why in Theorem \ref{phtr}(i) we are not able to prove
the full
range of (EG2), note that to prove that a configuration is good we must
show that
the color at $0$ cannot be affected by \textit{any} tampering of the
color record far
away from~$0$. In contrast, to prove that a configuration is bad it
suffices to
exhibit \textit{just two} tamperings that affect the color at $0$. In
essence, the
conditions on $p$ and $\eps$ in Theorem~\ref{phtr}(i) guarantee that
the random
walk has such a large drift that it moves away from the origin no
matter what the
color record is.

We close with (see Figure \ref{figbad}) the following conjecture.

\begin{conjecture}\label{phtrconj}
\textup{(EG2)} is true.
\end{conjecture}

Theorem \ref{phtr} is proved in Sections \ref{S2}--\ref{S5}: (i) in
Section \ref{S2},
(ii) in Section \ref{S4} and (iii) in Section \ref{S5}. It seems that
for $p \in
[\frac12,\frac45)$ and $\eps\in(0,1)$ the argument needed to
prove that all
configurations are bad is much more involved. In Section \ref{S3}, we
suggest an
approach to handle this problem, which will be pursued in future work.

The examples alluded to at the end of Section \ref{S1.1} typically
have both good and
bad configurations. On the other hand, we believe that our process $Y$
has all good
or all bad configurations, except at the point $(\frac12,0)$ and
possibly on the line
segment $\{\frac45\} \times[0,1)$. A simple example with such a
dichotomy, due to
Rob van den Berg, is the following. Let $X=(X_n)_{n\in\Z}$ be an
i.i.d. $\{0,1\}$-valued
process with the 1's having density $p\in(0,1)$. Let $Y_n=\ind\{
X_n=X_{n+1}\}$,
$n\in\Z$. Clearly, if $p=\frac12$, then $Y=(Y_n)_{n\in\Z}$ is
also i.i.d., and hence
all configurations are good. However, if $p \neq\frac12$, then it is
straightforward
to show that all configurations are bad. See \cite{LMV}, Proposition 3.3.


\section{\texorpdfstring{$\B=\varnothing$ for $p$ large and $\varepsilon$ small}{B = nothing for p large and epsilon small}}
\label{S2}

In this section, we prove Theorem~\ref{phtr}(i). The proof is based on
Lemmas \ref{l1}--\ref{l3}
in Section \ref{S2.1}, which are proved in Sections~\ref{S2.2}--\ref
{S2.4}, respectively. A key
ingredient of these lemmas is control of the \textit{cut times} for the
walk, that is, times at
which the past and the future of the walk have disjoint supports.
Throughout the paper, we abbreviate $I_m^n:=\{m,\ldots,n\}$ for
$m,n\in\N_0$ with $m\leq n$.


\subsection{\texorpdfstring{Proof of Theorem \protect\ref{phtr}(\textup{i}): Three lemmas}
{Proof of Theorem 1.3(i): Three lemmas}} \label{S2.1}

For $m,n\in\N$ with $m\leq n$, abbreviate
\[
S_m^n := (S_m,\ldots ,S_n) \quad\mbox{and}\quad Y_m^n := (Y_m,\ldots ,Y_n).
\]

The main ingredient in the proof of Theorem \ref{phtr}(i) will be an
estimate of the
number of cut times along $S_0^n$.

\begin{definition}\label{cuttimes}
For $n\in\N$, a time $k\in\N_0$ with $k\leq n-1$ is a cut time for
$S_0^n$ if
and only if
\[
S_0^k \cap S_{k+1}^n =\varnothing\quad\mbox{and}\quad S_k\geq0.
\]
\end{definition}

This definition takes into account only cut times corresponding to locations
on or to the right of the origin. Let $CT_n=CT_n(S_0^n)=CT_n(S_1^n)$
denote the set of cut times for $S_0^n$. Our first lemma reads as follows.

\begin{lemma}\label{l1}
For $k\in\N_0$, let $\cE_k\in\sigma(S_0^k,Y_0^k)$ be any event in
the $\sigma$-algebra
of the walk and the color record up to time $k$. Then
%
\begin{equation}
\label{Ekids}
\P(\cE_k \mid k\in CT_n, Y_1^n=y_1^n)
= \P(\cE_k \mid k\in CT_n,  Y_1^n=\bar y_1^n)
\end{equation}
for all $n \in\N$ with $n>k$ and all $y_1^n,\bar y_1^n$ such that
$y_1^k=\bar y_1^k$.\vspace*{3pt}
\end{lemma}

We next define
%
\begin{equation}
\label{fdef}
\hspace*{20pt}f(m) := \sup_{n\geq m} \max_{y_1^n}
\mathop{\mathop{\max}_{ {A\subseteq I_0^{m-1}} }}_{{|A| \geq {m}/{2}} }
 \P(CT_n\cap A =\varnothing\mid Y_1^n=y_1^n), \qquad m\in\N.
\end{equation}
Our second and third lemma read as follows.\vspace*{3pt}

\begin{lemma}\label{l2}
If $\lim_{m\rightarrow\infty} mf(m) =0$, then $\B=\varnothing$.\vspace*{3pt}
\end{lemma}

\begin{lemma}\label{l3}
$\limsup_{m\rightarrow\infty} \frac{1}{m} \log f(m) <0$ for
$(p,\eps)$ in a neighborhood of
$(1,0)$ containing the line segment $(p_*,1] \times\{0\}$.\vspace*{3pt}
\end{lemma}

Note that Lemma \ref{l3} yields the exponential decay of $m \mapsto
f(m)$, which is much more
than is needed in Lemma \ref{l2}. Note that Lemmas \ref{l2} and \ref
{l3} imply Theorem~\ref{phtr}(i).

Lemma \ref{l1} states that, conditioned on the occurrence of a cut
time at time $k$, the color
record after time $k$ does not affect the probability of any event that
is fully determined by
the walk and the color record up to time~$k$. Lemma \ref{l2} gives the
following sufficient
criterion for the nonexistence of bad configurations: for any set of
times up to time $m$ of
cardinality at least~$\frac{m}{2}$, the probability that the walk up
to time $n\ge m$ has no
cut times in this set, even when conditioned on the color record up to
time $n$, decays faster
than $\frac{1}{m}$ as $m\to\infty$, uniformly in $n$ and in the
color record that is being
conditioned on. Lemma \ref{l3} states that for $p$ and $\eps$ in the
appropriate range, the
above criterion is satisfied.

A key formula in the proof of Lemmas \ref{l1}--\ref{l3} is the
following. Let $R(s_1^n)$
denote the range of $s_1^n$ (i.e., the cardinality of its support), and
write $s_1^n \sim
y_1^n$ to denote that $s_1^n$ and $y_1^n$ are \textit{compatible} (i.e.,
there exists a coloring
of $\Z$ for which $s_1^n$ generates $y_1^n$). Below we abbreviate
$\P(S_1^n=s_1^n)$ by $\P(s_1^n)$.\vspace*{3pt}

\begin{proposition}\label{decomp}
For all $n\in\N$,
\[
\P(S_1^n=s_1^n,Y_1^n=y_1^n) = \P(s_1^n)  \bigl(\tfrac12
\bigr)^{R(s_1^n)}
\ind\{s_1^n \sim y_1^n\}.
\]
\vspace*{3pt}
\end{proposition}

The factor $(\frac12)^{R(s_1^n)}$ arises because if $s_1^n \sim
y_1^n$, then $y_1^n$ fixes
the coloring on the support of $s_1^n$.


\subsection{\texorpdfstring{Proof of Lemma \protect\ref{l1}}{Proof of Lemma 2.2}}
\label{S2.2}

Write $\P(\cE_k \mid k\in CT_n, Y_1^n=y_1^n) = N_k/D_k$ with (use
Proposition \ref{decomp})
\begin{eqnarray*}
N_k &:=& \sum_{x=0}^n \sum_{s_1^n} \ind\{s_k=x\} \ind\{ k\in
CT_n(s_1^n)\}
\P(s_1^n)  \biggl(\frac12 \biggr)^{R(s_1^n)} \ind\{s_1^n \sim
y_1^n\} \ind\{\cE_k\},\\
D_k &:=& \sum_{x=0}^n \sum_{s_1^n} \ind\{s_k=x\} \ind\{ k\in
CT_n(s_1^n)\}
\P(s_1^n)  \biggl(\frac12 \biggr)^{R(s_1^n)} \ind\{s_1^n \sim
y_1^n\}.
\end{eqnarray*}
Abbreviate $\{S_k^n>x\}$ for $\{S_l>x \ \forall k \leq l \leq n\}$,
etc. Note that if
$k\in CT_n(s_1^n)$, then we have
$\ind\{s_1^n \sim y_1^n\}= \ind\{s_1^k \sim y_1^k\} \ind\{s_{k+1}^n
\sim
y_{k+1}^n\}$ and $R(s_1^n) =R(s_1^k)+R(s_{k+1}^n)$. It follows that
\begin{eqnarray*}
N_k &=& \sum_{x=0}^n \sum_{s_1^k} \ind\{s_k=x\} \ind\{s_1^k\leq x\}
\P(s_1^k)  \biggl(\frac12 \biggr)^{R(s_1^k)} \ind\{s_1^k \sim
y_1^k\}\ind\{\cE_k\} \\
&&\hphantom{\sum_{x=0}^n}{} \times\sum_{s_{k+1}^n} \ind\{s_{k+1}^n>x\} \P
(s_{k+1}^n\mid S_k=x)
 \biggl(\frac12 \biggr)^{R(s_{k+1}^n)} \ind\{s_{k+1}^n \sim
y_{k+1}^n\}\\
&=& C_{k,n}(y_{k+1}^n) \sum_{x=0}^n \sum_{s_1^k} \ind\{s_k=x\} \ind
\{ s_1^k\leq x\}
\P(s_1^k)  \biggl(\frac12 \biggr)^{R(s_1^k)} \ind\{s_1^k \sim
y_1^k\} \ind\{\cE_k\}
\end{eqnarray*}
with (shift $S_k$ back to the origin)
\[
C_{k,n}(y_{k+1}^n) := \biggl [ \sum_{s_{1}^{n-k}} \ind\{s_{1}^{n-k}>0\}
\P(s_{1}^{n-k})  \biggl(\frac12 \biggr)^{R(s_1^{n-k})}
\ind\{s_{1}^{n-k} \sim y_{k+1}^n\}  \biggr].
\]
Likewise, we have
\[
D_k = C_{k,n}(y_{k+1}^n) \sum_{x=0}^n \sum_{s_1^k} \ind\{s_k=x\}
\ind\{ s_1^k\leq x\}
\P(s_1^k) \biggl (\frac12 \biggr)^{R(s_1^k)} \ind\{s_1^k \sim
y_1^k\}.
\]
The common factor $C_{k,n}(y_{k+1}^n)$ cancels out and so $N_k/D_k$
only depends on~$y_1^k$.
Therefore, as long as $y_1^k=\bar y_1^k$, we have the equality in (\ref
{Ekids}).\vspace*{12pt}


\subsection{\texorpdfstring{Proof of Lemma \protect\ref{l2}}{Proof of Lemma 2.3}}
\label{S2.3}

Since $f(m)\leq\frac12$ for all large $m$, we will assume that
all the values of $m$ arising in the proof below satisfy this.

For $n\in\N$ and $y_1^n$ and $\bar y_1^n$, define
\[
\Delta^n(y_1^n,\bar y_1^n) := \P(Y_0=B \mid Y_1^n=y_1^n) - \P(Y_0=B
\mid Y_1^n =\bar y_1^n).
\]
We will show that if $\lim_{n\rightarrow\infty} mf(m)=0$, then
%
\begin{equation}
\label{mainest}
\lim_{m\to\infty} \sup_{n\geq m}
\mathop{\mathop{\max}_{ {y_1^n,\bar y_1^n} }}_{{y_1^{m-1}=\bar y_1^{m-1}} }
|\Delta^n(y_1^n,\bar y_1^n) |= 0,
\end{equation}
and hence $\B=\varnothing$ by Definition \ref{badconf}.

In what follows, we
%
\begin{equation}
\label{fixY}
\mbox{fix }   m,n\in\N  \mbox{ with }   m\leq n   \mbox{ and }
y_1^n,\bar y_1^n   \mbox{ with }   y_1^{m-1}=\bar y_1^{m-1}
\end{equation}
and abbreviate $\Delta=\Delta^n(y_1^n,\bar y_1^n)$. Define
\begin{eqnarray*}
A &\hspace*{3pt}=& A_m^n (y_1^n,\bar y_1^n)\\
&:=&  \{k\in I_0^m\dvtx
\P(k\in CT_n \mid Y_1^n=y_1^n) - \P(k\in CT_n \mid Y_1^n =\bar y_1^n)
\geq- 2f(m) \}.
\end{eqnarray*}
Using Lemma \ref{l1}, we will show that
%
\begin{equation}
\label{stringimpl}
|A| \geq\frac{m}{2}
\end{equation}
and
%
\begin{equation}
\label{stringimplAGAIN}
|\Delta| \leq2f(m)(m+1).
\end{equation}
The argument we will give works for any choice of $y_1^n$ and $\bar y_1^n$
subject to (\ref{fixY}) (with the corresponding $A$ and $\Delta$).
Together with $\lim_{m\rightarrow\infty} mf(m) =0$, (\ref{stringimplAGAIN})
will prove Lemma \ref{l2}.\vspace*{12pt}

\subsubsection{\texorpdfstring{Proof of (\protect\ref{stringimpl})}{Proof of (2.5)}}
\label{S2.3.2}

Write $B:=I_0^{m-1} \setminus A = \{b_1,\ldots ,b_{m-|A|}\}$. We will
show that
$f(m)\leq\frac{1}{2}$ and $|B|>\frac{m}{2}$ are incompatible.
Indeed, by the definition
of $A$, we have
\begin{eqnarray}
\P(b_i\in CT_n \mid Y_1^n=y_1^n) - \P(b_i\in CT_n \mid Y_1^n= \bar y_1^n)
< -2f(m),\nonumber\\
 \eqntext{i=1,\ldots ,m-|A|.}
\end{eqnarray}
Define $B_i:=\{b_1,\ldots ,b_i\}$, $i=1,\ldots ,m-|A|$, with the
convention that $B_0=\varnothing$.
Estimate, writing $\mathit{FCT}_n(B)$ to denote the first cut time for $S_0^n$
in $B$,
\begin{eqnarray*}
&&\P(CT_n\cap B \neq\varnothing\mid Y_1^n=y_1^n) - \P(CT_n\cap B \neq
\varnothing
\mid Y_1^n =\bar y_1^n)\\
&& \qquad  = \sum_{i=1}^{m-|A|}
\bigl [\P\bigl(\mathit{FCT}_n(B)=b_i\mid Y_1^n=y_1^n\bigr) - \P\bigl(\mathit{FCT}_n(B)=b_i \mid
Y_1^n=\bar y_1^n\bigr) \bigr]\\
&& \qquad  = \sum_{i=1}^{m-|A|}
\P (CT_n\cap B_{i-1} = \varnothing\mid b_i\in CT_n ,
Y_1^n=y_1^n )\\
&&\hphantom{\sum_{i=1}^{m-|A|}} \qquad  \quad {}   \times [ \P(b_i\in CT_n \mid Y_1^n=y_1^n)
- \P(b_i\in CT_n \mid Y_1^n =\bar y_1^n)  ]\\
&& \qquad  < -2f(m) \sum_{i=1}^{m-|A|} \P (CT_n\cap B_{i-1}
= \varnothing\mid b_i\in CT_n, Y_1^n=y_1^n )\\
&& \qquad  \leq-2f(m) \sum_{i=1}^{m-|A|} \P (b_i\in CT_n ,
CT_n\cap B_{i-1}
= \varnothing\mid Y_1^n=y_1^n )\\
&& \qquad  = - 2f(m)  [ 1-\P(B\cap CT_n = \varnothing\mid
Y_1^n=y_1^n) ]  ,
\end{eqnarray*}
where in the third line we have used Lemma \ref{l1}. This inequality
can be rewritten as
\[
2f(m) < \P(CT_n\cap B = \varnothing\mid Y_1^n=y_1^n) \bigl(1+2f(m)\bigr)
- \P(CT_n\cap B = \varnothing\mid Y_1^n = \bar y_1^n).
\]
By (\ref{fdef}), the right-hand side is at most $f(m)(1+2f(m))$ when
$|B|>\frac{m}{2}$,
which gives a contradiction because $f(m)\leq\frac{1}{2}$.

\subsubsection{\texorpdfstring{Proof of (\protect\ref{stringimplAGAIN})}{Proof of (2.6)}}
\label{S2.3.1}

Write
\[
\tilde{\Delta} := \P (Y_0=B, CT_n\cap A \neq\varnothing\mid
Y_1^n=y_1^n )
- \P (Y_0=B, CT_n\cap A \neq\varnothing\mid Y_1^n =\bar y_1^n ).
\]
Using (\ref{fdef}) in combination with (\ref{stringimpl}),
we may estimate
\[
\Delta\leq
\tilde{\Delta} + f(m).
\]
Let $A=\{a_1,\ldots,a_{|A|}\}$ denote the elements of $A$ in
increasing order,
and define $A_i:=\{a_1,\ldots,a_i\}$, $i=1,\ldots,|A|$, with the convention
that $A_0=\varnothing$. Then, using Lemma \ref{l1}, we have
\begin{eqnarray*}
\tilde{\Delta} &=& \sum_{i=1}^{|A|}  \bigl[\P \bigl(Y_0=B,
\mathit{FCT}_n(A)=a_i \mid Y_1^n=y_1^n \bigr)
\\[-1pt]
&&\hphantom{\sum_{i=1}^{|A|}  \bigl[}{}- \P\bigl (Y_0=B, \mathit{FCT}_n(A)=a_i \mid Y_1^n =\bar y_1^n \bigr) \bigr]\\[-1pt]
&=& \sum_{i=1}^{|A|} [ \P (Y_0=B,  CT_n\cap A_{i-1}
= \varnothing\mid a_i \in CT_n, Y_1^n=y_1^n )
  \P(a_i\in CT_n \mid Y_1^n=y_1^n)\\[-1pt]
&&\hphantom{\sum_{i=1}^{|A|} [ }{}  - \P (Y_0=B,  CT_n\cap A_{i-1}
= \varnothing\mid a_i\in CT_n,  Y_1^n=\bar y_1^n )\\
&&\hspace*{188pt}\hphantom{\sum_{i=1}^{|A|} [ - \P }{}\times
  \P(a_i\in CT_n \mid Y_1^n= \bar y_1^n)  ]\\[-3.5pt]
&=& \sum_{i=1}^{|A|} \P (Y_0=B,  CT_n\cap A_{i-1}
= \varnothing\mid a_i\in CT_n, Y_1^n=y_1^n )  D_i,
\end{eqnarray*}
where\vspace*{-1.5pt}
\[
D_i := \P(a_i\in CT_n \mid Y_1^n=y_1^n) - \P(a_i\in CT_n \mid Y_1^n
=\bar y_1^n).
\]
In the third line, we have used the fact that $\{CT_n\cap A_{i-1}
=\varnothing\} = \{A_{i-1}
\cap CT_{a_i} =\varnothing\} \in\sigma(S_0^{a_i},Y_0^{a_i})$ (the
$\sigma$-algebra generated
by $S_0^{a_i},Y_0^{a_i}$) on the event $\{a_i\in CT_n\}$, so that Lemma
\ref{l1} applies.
The definition of the set $A$ implies that $D_i \geq-2f(m)$ for all
$i$. Hence, by using
Lemma \ref{l1} once more, we obtain\vspace*{-1pt}
\begin{eqnarray*}
\tilde{\Delta}  &\leq&\sum_{i=1}^{|A|} \ind\{D_i\geq0\}
\P (Y_0=B, CT_n\cap A_{i-1} = \varnothing\mid a_i \in CT_n,
Y_1^n=y_1^n ) D_i\\[-1pt]
   &\leq&\sum_{i=1}^{|A|} \ind\{D_i\geq0\}
\P (CT_n\cap A_{i-1} = \varnothing\mid a_i\in CT_n,
Y_1^n=y_1^n )  D_i\\
  &=& \sum_{i=1}^{|A|} \P (CT_n\cap A_{i-1} = \varnothing\mid a_i \in
CT_n, Y_1^n=y_1^n ) D_i\\
&&  {}  + \sum_{i=1}^{|A|} \ind\{D_i<0\}
\P (CT_n\cap A_{i-1} = \varnothing\mid a_i \in CT_n,
Y_1^n=y_1^n ) (-D_i)\\
  &\leq&\sum_{i=1}^{|A|}  [\P (a_i\in CT_n, CT_n\cap A_{i-1}
= \varnothing\mid Y_1^n=y_1^n )\\
&&\hphantom{\sum_{i=1}^{|A|}  [}{}
- \P (a_i\in CT_n, CT_n\cap A_{i-1} = \varnothing\mid Y_1^n= \bar
y_1^n ) ]   + 2f(m)|A|\\
  &=& \P(CT_n\cap A \neq\varnothing\mid Y_1^n=y_1^n) - \P(CT_n\cap A
\neq\varnothing
\mid Y_1^n= \bar y_1^n) + 2f(m)|A|\\
   &\leq& f(m) + 2f(m)m.
\end{eqnarray*}
Thus, we find that $\Delta\leq2f(m)(m+1)$, where the upper bound does
not depend on the
choice of configurations made in (\ref{fixY}). Exchanging $y_1^n$ and~$\bar y_1^n$,
we obtain the same bound for $|\Delta|$.
Hence, we have proved (\ref{stringimplAGAIN}).


\subsection{\texorpdfstring{Proof of Lemma \protect\ref{l3}}{Proof of Lemma 2.4}}
\label{S2.4}

For simplicity, we will only consider $m$-values that are a multiple of~$6$.
The proof is easily adapted to intermediate $m$-va\-lues.

We first state the following fairly straightforward lemma, where we
note that
$\{S_m^n > \frac{2m}{3}\} = \{S_l > \frac{2m}{3} \ \forall m
\leq l\leq n\}$.

\begin{lemma}\label{easy}
For $m,n\in\N$ with $m\leq n$,
%
\begin{equation}
\label{TCspeed}
\biggl\{|CT_n\cap I_0^{m-1}| \leq\frac{m}{2}\biggr\}
\subseteq\biggl\{S_m^n > \frac{2m}{3}\biggr\}^c.
\end{equation}
\end{lemma}

\begin{pf}
Note that each cut time $k$ corresponds to a cut point $S_k$,
and so the set $CT_n\cap I_0^{m-1}$ of cut times corresponds to a set
$CP_n(m)$ of cut
points. On the event $\{S_m^n>\frac{2m}{3}\}$, the interval
$I_0^{ {2m}/{3}}$
is fully covered by $S_0^{m-1}$. For each $x\in I_0^{ {2m}/{3}}$, we
look at the steps of the random walk entering or exiting~$x$ from the right:
\begin{itemize}
\item
If $x\in CP_n(m)$, then during the time interval $I_0^{n-1}$
there is at least one step exiting $x$ to the right.
\item
If $x\notin CP_n(m)$, then during the time interval $I_0^{n-1}$
there are at least two
steps exiting $x$ to the right and one step entering $x$ from the right
(since there must be
a return to $x$ from the right).
\end{itemize}
Since each step refers to a single point $x$ only, and $S_0^{m-1}$ goes
along at most~$m$ edges (and exactly $m$ edges when $\eps=0$), we get that
%
\[
m \geq |CP_n(n)\cap I_0^{ {2m}/{3}} | +
3  |I_0^{ {2m}/{3}}
\bk CP_n(n) |
= 3\biggl(\frac{2m}{3} + 1\biggr) - 2  |CP_n(n)\cap I_0^{ {2m}/{3}} |.
\]
Hence, $|CP_n(n)\cap I_0^{ {2m}/{3}}|>\frac{m}{2}$. Still on the event
$\{S_m^n>\frac{2m}{3}\}$, the cut times corresponding to $CP_n(n)\cap
I_0^{ {2m}/{3}}$ occur before time $m-1$, and so
\[
 |CT_n \cap I_0^{m-1} | \geq |CP_n(n)\cap I_0^{
{2m}/{3}} |.
\]
Hence, $|CT_n\cap I_0^{m-1}|>\frac{m}{2}$, and so (\ref{TCspeed}) is proved.
\end{pf}

For $A\subseteq I_0^{m-1}$ such that $|A|\geq\frac{m}{2}$, we have
\[
\{CT_n\cap A =\varnothing\} \subseteq \biggl\{ |CT_n\cap I_0^{m-1}
 | \leq\frac{m}{2} \biggr\}.
\]
Therefore, by (\ref{TCspeed}),
%
\begin{eqnarray}
\label{lastvisit}
\hspace*{20pt}\{CT_n\cap A =\varnothing\}
&\subseteq& \biggl\{\exists k\dvtx  m\leq k \leq n-1,   S_k = \frac{2m}{3},
  S_{k+1}^n >\frac{2m}{3} \biggr\} \nonumber
  \\[-8pt]
  \\[-8pt]
  \hspace*{20pt} &&{}\cup\biggl\{S_n \leq\frac{2m}{3}\biggr\}.
\nonumber
\end{eqnarray}


\vspace*{12pt}
\subsubsection{\texorpdfstring{Estimate of the probabilities of the events in (\protect\ref{lastvisit})}
{Estimate of the probabilities of the events in (2.8)}}
\label{S2.4.1}

In this subsection, we obtain upper bounds on the probabilities
of the two events on the right-hand side of (\ref{lastvisit}) when
conditioned on $Y_1^n$. The upper bounds will appear in
(\ref{1stterm}) and (\ref{2ndterm}) below. In Section \ref{S2.4.2},
we use these estimates to finish the proof of Lemma~\ref{l3}.

Write
%
\begin{eqnarray}
\label{decomp2}
&&\P \biggl(\exists k\dvtx  m\leq k \leq n-1,   S_k
= \frac{2m}{3}, S_{k+1}^n >\frac{2m}{3}\Bigm| Y_1^n=y_1^n \biggr)\nonumber
\\[-8pt]
\\[-8pt]
&& \qquad  = \sum_{k=m}^{n-1} \P \biggl(S_{k}=\frac{2m}{3},
S_{k+1}^n >\frac{2m}{3} \Bigm| Y_1^n=y_1^n \biggr)
= \sum_{k=m}^{n-1} \frac{N_k}{D_k},
\nonumber
\end{eqnarray}
with (recall Proposition \ref{decomp})
\begin{eqnarray*}
N_k &:=& N_k(y_1^n) = \sum_{s_1^n} \ind\biggl\{s_{k}=\frac{2m}{3}\biggr\} \ind
\biggl\{s_{k+1}^n
>\frac{2m}{3}\biggr\} \P(s_1^n)
\biggl (\frac12 \biggr)^{R(s_1^n)} \ind\{s_1^n \sim y_1^n\},\\
D_k &:=& D_k(y_1^n) = \sum_{s_1^n} \P(s_1^n)
 \biggl(\frac12 \biggr)^{R(s_1^n)} \ind\{s_1^n \sim y_1^n\}.
\end{eqnarray*}
Estimate
\begin{eqnarray*}
&&N_k \leq\sum_{s_1^k} \ind\biggl\{s_k=\frac{2m}{3}\biggr\} \P(s_1^k) \ind\{
s_1^k \sim y_1^k\}\\
&& \qquad {}  \times\sum_{s_{k+1}^n} \ind\biggl\{s_{k+1}^n>\frac
{2m}{3}\biggr\}
\P\biggl(s_{k+1}^n \Bigm| S_k=\frac{2m}{3}\biggr)
\biggl (\frac12 \biggr)^{R(s_{k+1}^n)} \ind\{s_{k+1}^n \sim
y_{k+1}^n\}.
\end{eqnarray*}
Here, the bound arises by noting that $\ind\{s_1^n \sim y_1^n\} \leq
\ind\{s_1^k
\sim y_1^k\}\ind\{s_{k+1}^n \sim y_{k+1}^n\}$ and estimating $R(s_1^n)
\geq
R(s_{k+1}^n)$. Thus, shifting $S_k$ back to the origin, we get
%
\begin{equation}
\label{Nkub}
N_k \leq\P \biggl(S_k=\frac{2m}{3},S_1^k \sim y_1^k \biggr) C_{k,n}(y_{k+1}^n)
\end{equation}
with
\[
C_{k,n}(y_{k+1}^n) = \sum_{s_1^{n-k}}
\ind\{s_1^{n-k}>0\} \P(s_1^{n-k}) \biggl (\frac12
\biggr)^{R(s_1^{n-k})}
\ind\{s_1^{n-k} \sim y_{k+1}^n\}.
\]
Next, estimate
\begin{eqnarray*}
D_k &\geq&\sum_{s_1^k} \ind\{s_1^k \leq s_k\}
\P(s_1^k)  \biggl(\frac12 \biggr)^{R(s_1^k)} \ind\{s_1^k \sim
y_1^k\}\\
&&{} \times\sum_{s_{k+1}^n} \ind\{s_{k+1}^n >s_k\} \P(s_{k+1}^n
\mid S_k=s_k)
 \biggl(\frac12 \biggr)^{R(s_{k+1}^n)} \ind\{s_{k+1}^n \sim
y_{k+1}^n\}.
\end{eqnarray*}
Here, the bound arises by restricting $S_1^n$ to the event
\[
\{k\in CT_n\} = \{S_1^k \leq S_k\} \cap\{S_{k+1}^n >S_k\},
\]
noting that $\ind\{S_1^n \sim y_1^n\}=\ind\{S_1^k \sim y_1^k\}\ind\{S_{k+1}^n
\sim y_{k+1}^n\}$ on this event, and inserting $R(s_1^n) = R(s_1^k) +
R(s_{k+1}^n)$. Thus, shifting $S_k$ back to the origin, we get
%
\begin{equation}
\label{Dklb}
D_k \geq\E\bigl (\bigl (\tfrac12 \bigr)^{R(S_1^k)} \ind\{S_1^k
\leq S_k\}
\ind\{S_1^k \sim y_1^k\} \bigr) C_{k,n}(y_{k+1}^n).
\end{equation}
Combining the upper bound on $N_k$ in (\ref{Nkub}) with the lower
bound on $D_k$ in~(\ref{Dklb}),
and canceling out the common factor $C_{k,n}(y_{k+1}^n)$, we arrive at
%
\begin{eqnarray}
\label{NDk}
 && \P \biggl(S_k=\frac{2m}{3},   S_{k+1}^n >\frac{2m}{3} \Bigm|
Y_1^n=y_1^n \biggr)\nonumber
\\[-8pt]
\\[-8pt]
&&\qquad\leq\frac{\P(S_k= {2m}/{3},S_1^k \sim y_1^k)}{\E (
( 1/2 )^{R(S_1^k)}
\ind\{S_1^k \leq S_k\} \ind\{S_1^k \sim y_1^k\} )}.
\nonumber
\end{eqnarray}
Note that this bound is uniform in $n$.

The numerator of (\ref{NDk}) is bounded from above by $\P(S_k=\frac
{2m}{3})$, while the
denominator of (\ref{NDk}) is bounded from below by $(\frac12)^k\P(S_k=k)
=(\frac{p(1-\eps)}{2})^k$, where we note that $S_1^k \sim y_1^k$ for
all $y_1^k$ on the
event $\{S_k=k\}$. Hence, by (\ref{decomp2}), we have
%
\begin{eqnarray}
\label{1stterm}
&&\P \biggl(\exists k\dvtx  m\leq k \le n-1 ,   S_k = \frac{2m}{3},
S_{k+1}^n >\frac{2m}{3} \Bigm| Y_1^n=y_1^n \biggr)
\nonumber
\\[-8pt]
\\[-8pt]&& \qquad \leq\sum_{k=m}^{n-1} \frac{\P(S_k =  {2m}/{3})}{(
{p(1-\eps)}/{2})^k}.
\nonumber
\end{eqnarray}
The bound in (\ref{1stterm}) controls the first term in the right-hand
side of (\ref{lastvisit}).

Since $\P(Y_1^n=y_1^n) \geq\P(Y_1^n=y_1^n,S_n=n) = (\frac{p(1-\eps
)}{2})^n$, we have
%
\begin{equation}
\label{2ndterm}
\hspace*{20pt}\P\biggl(S_n\leq\frac{2m}{3} \Bigm| Y_1^n=y_1^n\biggr) \leq\frac{\P(S_n \leq
 {2m}/{3})}
{( {p(1-\eps)}/{2})^n}
\leq C  \frac{\P(S_n =  {2m}/{3})}{( {p(1-\eps)}/{2})^n},
\end{equation}
provided $n$ is even (which is necessary when $\eps=0$ because we have
assumed that
$\frac{2m}{3}$ is even). Here, the constant $C =C(p,\varepsilon) \in
(1,\infty)$ comes from
an elementary large deviation estimate, for which we must assume that
%
\begin{equation}
\label{driftcond}
(2p-1)(1-\eps) > \tfrac23.
\end{equation}
The bound in (\ref{2ndterm}) controls
the second term in the right-hand side of (\ref{lastvisit}).

\subsubsection{Completion of the proof}
\label{S2.4.2}

In this section, we finally complete the proof of Lemma \ref{l3}.

Combining (\ref{1stterm})--(\ref{2ndterm}) and recalling (\ref{fdef})
and (\ref{lastvisit}), we obtain the estimate
%
\begin{equation}
\label{ACTbdalt}
f(m) \leq(C+1) \sum_{k= {m}/{2}}^{\infty}
\frac{\P(S_{2k} =  {2m}/{3})}{( {p(1-\eps)}/{2})^{2k}}.
\end{equation}
Since there exists a $C'=C'(p,\eps)\in(1,\infty)$ such that, for
$k\geq\frac12 m$,
\[
\P\biggl(S_{2k} = \frac{2m}{3}\biggr)\le C'\P\biggl(S_{2k} = \frac{4k}{3}\biggr),
\]
we see that $\limsup_{m\to\infty}\frac{1}{m}\log f(m)<0$ as soon as
%
\begin{equation}
\label{ACTbd}
\limsup_{m\to\infty} \frac{1}{m} \log\P\biggl(S_m = \frac{2m}{3}\biggr)
< \log\biggl (\frac{p(1-\eps)}{2} \biggr).
\end{equation}

Note that (\ref{driftcond}) holds for $(p,\eps)$ in a neighborhood of $(1,0)$
containing the line segment $(p_*,1] \times\{0\}$.

By Cramer's theorem of large deviation theory (see, e.g., \cite{H00},
Chapter~I), the left-hand
side of (\ref{ACTbd}) equals $-I(p,\eps)$ with
%
\begin{equation}
\label{Ipespdef}
I(p,\eps) := \sup_{\lambda\in\R} \biggl [\frac23\lambda- \log
M(\lambda;p,\eps) \biggr],
\end{equation}
where
%
\begin{equation}
\label{Mlampepsdef}
M(\lambda;p,\eps) := \eps+ p(1-\eps)e^\lambda+ (1-p)(1-\eps
)e^{-\lambda}
\end{equation}
is the moment-generating function of the increments of $S$. Due to the
strict convexity of
$\lambda\mapsto\log M(\lambda;p,\eps)$, the supremum is attained at
the unique $\bar\lambda$
solving the equation
%
\begin{equation}
\label{barlam}
\frac23 = \frac{ (\partial/\partial\lambda) M(\lambda;p,\eps)
}{ M(\lambda;p,\eps) },
\end{equation}
where we note that $\bar\lambda< 0$ because of (\ref{driftcond}).
For the special case where
$\eps= 0$, an easy calculation gives
\[
\bar\lambda= \frac12 \log \biggl(\frac{5(1-p)}{p} \biggr),
\]
implying that $I(p,0) = \log C(p)$ with $C(p) =
[5/6p]^{5/6}[1/6(1-p)]^{1/6}$. Hence, the
inequality in (\ref{ACTbd}) reduces to $C(p)> 2/p$, which is
equivalent to $p>p^*$ with
$p^*=1/(1+5^5 12^{-6})$. The same formulas (\ref{Ipespdef})--(\ref
{barlam}) show that
(\ref{ACTbd}) holds in a neighborhood of $(1,0)$.\vspace*{12pt}


\section{\texorpdfstring{$\B=\Omega$ for $p \in(\frac12,\frac45)$ and $\varepsilon=0$}
{B = Omega for p in (1/2,4/5) and epsilon = 0}}
\label{S4}

Throughout the remainder of this paper [with the sole
exceptions of Section \ref{S5.1}
and the claim of independence immediately prior to (\ref{lim1})],
we use $Y_1^\infty$, $\bar{Y}_1^\infty$
and $\tilde{Y}_1^\infty$ to represent specific sequences rather
than random sequences. This abuse of notation will nowhere cause
harm.

In this section, we prove Theorem \ref{phtr}(ii). The proof is based
on the following
observations valid for a random walk that cannot pause ($\eps=0$).

\begin{longlist}[(II)]
\item[(I)]
On a color record of the type $[\mathit{WWBB}]^M$, $M\in\N$, the walk cannot
turn. Indeed, a turn
forces the same color to appear in the color record two units of time apart.
\item[(II)]
Any color record $Y_1^{m-1}$ up to time $m\in\N$ can be seen in a
unique way along a
stretch of coloring of the type $[\mathit{WWBB}]^M$ with $M \geq m$. Indeed, on
such a
stretch each site has a $W$-neighbor and a $B$-neighbor, so once the
starting or
ending point of the walk is fixed it is fully determined by $Y_1^{m-1}$.
\end{longlist}

We prove Theorem \ref{phtr}(ii) by showing the following claim:
\begin{itemize}
\item
For any $Y_1^\infty$, $p \in(\frac12,\frac45)$ and $m\in\N$, we
can find $\bar Y_m^\infty$
and $\widetilde Y_m^\infty$ such that
%
\begin{equation}
\label{limbartil}
\hspace*{20pt}\lim_{n\to\infty}  \bigl|\P (C_0=W \mid Y_1^{m-1} \vee\bar
Y_m^n )
- \P (C_0=W \mid Y_1^{m-1} \vee\widetilde Y_m^n ) \bigr| = 2p-1,
\end{equation}
\end{itemize}
where $\vee$ denotes the concatenation operation. In view of Definition \ref{badconf}, this claim will imply that
$Y_1^\infty$ is bad.


\begin{pf}
Fix $m\in\N$.

1. We begin with the choice of $\bar Y_m^n$. For $L\in\N$, let
%
\begin{eqnarray}
\label{Yblocks}
\bar Y_m^n &:=& [\mathit{WWBB}]^m \mathit{WBB} [\mathit{WWBB}]^{2m} \mathit{WBB}
[\mathit{WWBB}]^{2m+1}\nonumber
\\[-8pt]
\\[-8pt]
&&{} \cdots \mathit{WBB} [\mathit{\mathit{WWBB}}]^{2m+L-1} \mathit{WBB} [\mathit{WWBB}]^{2m+L}.
\nonumber
\end{eqnarray}
The interest in this color record relies on three facts:
\begin{longlist}[(3)]
\item[(1)]
For $l=0,\ldots,L$, on the color record $[\mathit{WWBB}]^{2m+l}$ the walk
cannot turn
[see (I) above].
\item[(2)]
On $\bar Y_m^n$, the isolated $W$'s at the beginning of the $\mathit{WBB}$'s
play the
role of \textit{pivots}, since the walk can only turn there as is easily checked.
We call $W_0$ the pivot $W$ seen at time $5m$ (this is the first
pivot) and $W_l$, $l=1,\ldots,L$, the
subsequent pivots seen at times
\[
t(l) := 5m + \sum_{j=0}^{l-1} [3+4(2m+j)]
= k(2k+8m+1) + 5m, \qquad k=1,\ldots,L.
\]
\item[(3)]
Since the length of the color record $[\mathit{WWBB}]^{2m+l}$ increases with
$l$, if the walk does not turn on pivot $W_l$, then it cannot turn on
any later pivot. Indeed, going straight through $W_l$ means that the
coloring has an isolated~$W$ surrounded by two $B$'s, and this color
stretch is impossible to cross at any later time with any color record
of the
type $[\mathit{WWBB}]^M$, $M\in\N$.
\end{longlist}
The first color record $[\mathit{WWBB}]^m$ serves to prevent $W_0$ from being
in the coloring seen by the walk up to time $m-1$, because the walk
cannot turn between time $m$ and time $5m$ [see (I) above]. The total
time is
\[
n=n(L)=L(2L+8m+5)+13m+2.
\]

The above three facts imply that the behavior of the walk from time $m$
to time~$n$ (i.e.,
the increments $X_{m+1},\ldots,X_n$), leading to $\bar Y_m^n$ as its
color record, can be
characterized by the first pivot $W_l$, if any, where the walk makes no
turn. There are
$L+2$ possibilities, including the ones where there is a turn at every
pivot or at no pivot.
This characterization is up to a~2-fold symmetry in the direction of
the last step of the
walk, which can be either upwards or downwards (this is the same
symmetry as $X \to-X$).
Note that, except for the case where the walk makes no turn from time~$m$ to time $n$, the
behavior of the walk from time $1$ to time $m$ (i.e., the increments
$X_2,\ldots,X_m$) is
fully determined (up to the 2-fold symmetry) by $\bar Y_1^n$ [see (II)
above]. This is
because $l \mapsto t(l+1)-t(l)$ is increasing, so that $t(l+1)-t(l)
\geq t(1)-t(0)=3+8m>5m$.

Our goal will be to prove that for large $L$ the walk, conditioned on
$Y_1^{m-1} \vee
\bar Y_m^n$, with a high probability turns on every pivot and ends by
moving upwards.
To that end, we define the following events for the walk up to time~$n$:
\begin{itemize}
\item
$LT_l := \{$the walk turns on pivots $W_0,W_1,\ldots,W_l$ and does not
turn on pivots
$W_{l+1},\ldots,W_L\}$ (``last turn on $l$''), $l=0,\ldots,L$.
\item
$NT := \{\mbox{the walk does not turn on any pivot}\}$ (``no turn'').
\item
$EU:=\{S_n=S_{n-1}+1\}$ (``end upwards'').
\item
$ED:=\{S_n=S_{n-1}-1\}$ (``end downwards'').
\end{itemize}
Using these events, we may write
%
\begin{eqnarray}
\label{1split}
\hspace*{20pt}1 &=& \P (NT,EU \mid Y_1^{m-1} \vee\bar Y_m^n )
+ \P (NT,ED \mid Y_1^{m-1} \vee\bar Y_m^n )\nonumber
\\[-8pt]
\\[-8pt]
\hspace*{20pt}&&{}+ \sum_{l=0}^L  [\P (LT_l,EU \mid Y_1^{m-1} \vee
\bar Y_m^n )
+ \P (LT_l,ED \mid Y_1^{m-1} \vee\bar Y_m^n ) ].
\nonumber
\end{eqnarray}

Now, on the event $LT_l$, the length of the coloring seen by the walk
from time $1$ to time $n$ is
\[
n-t(l)+1 = \sum_{j=l}^L [3+4(2m+j)] = (L-l+1)(2L+2l+8m+3).
\]
Only two walks from time $m$ to time $n$ are in $LT_l$ and these are
reflections of each
other (one in $EU$ and one in $ED$). For either of these two walks, we
have that
$|S_{t(l)}-S_{t(0)}|=u(l)$, where
\[
u(l) := \sum_{j=1}^l (-1)^{l-j}[t(j)-t(j-1)]
= (1+8m) 1\{l\mbox{ odd}\} + 2l.
\]
It is easily checked that any walk in $EU\cap LT_l$ ends a distance at
least $2v(l,L)$
above any walk in $ED\cap LT_l$, with (see Figure \ref{figzigzag})
\begin{eqnarray*}
v(l,L) &:=& n(L)-t(l)-u(l)+(-1)^{l-1}4m - m\\
&\hspace*{3pt}=& (L-l)(2L+2l+8m+5)+2l+3m+2-1\{l\mbox{ odd}\}.
\end{eqnarray*}
Hence we have, using the fact that all walks in $LT_l$ visit the same
number of colors,
%
\begin{equation}
\label{ineqkk0}
 \hspace*{20pt}\quad \P (LT_l,ED \mid Y_1^{m-1} \vee\bar Y_m^n )
\leq \biggl(\frac{1-p}{p}\biggr )^{v(l,L)}
\P (LT_l,EU \mid Y_1^{m-1} \vee\bar Y_m^n ).
\end{equation}
Since $(1-p)/p<1$ (because $p>\frac12$) and $\lim_{L\to\infty} \inf
_{0\leq l \leq L}
v(l,L)=\infty$, it follows that for $L$ large the probability of
$LT_l\cap ED$ is
negligible with respect to the probability of $LT_l\cap EU$ uniformly
in $l$.


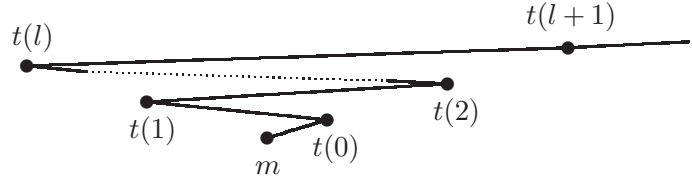
\begin{figure}
\begin{center}
\vspace*{-1cm}
\setlength{\unitlength}{0.4cm}
\hspace*{1.6cm}
\begin{picture}(10,10)(0,-2)
\put(0,0){\circle*{0.4}}
\put(2,.6){\circle*{0.4}}
\put(-4,1.2){\circle*{0.4}}
\put(6,1.8){\circle*{0.4}}
\put(-8,2.4){\circle*{0.4}}
\put(10,3){\circle*{0.4}}
\put(-.4,-1.2){$m$}
\put(1.5,-.6){$t(0)$}
\put(-4.6,0){$t(1)$}
\put(5.5,.6){$t(2)$}
\put(-8.5,3.1){$t(l)$}
\put(8.6,3.8){$t(l+1)$}
{\thicklines
\qbezier(0,0)(1,.3)(2,.6)
\qbezier(2,.6)(-1,.9)(-4,1.2)
\qbezier(-4,1.2)(1,1.5)(6,1.8)
\qbezier(6,1.8)(5,1.85)(4,1.9)
\qbezier(-6,2.2)(-7,2.3)(-8,2.4)
\qbezier(-8,2.4)(1,2.7)(10,3)
\qbezier(10,3)(12,3.1)(14,3.2)
}
\qbezier[50](4,1.9)(-1,2.05)(-6,2.2)
\end{picture}
\end{center}
\caption{A walk in $LT_l \cap EU$. The last turn occurs
at time
$t(l)$. Depending on the parity of $l$, the walk between time $m$ and time
$t(l)$ starts its zigzag motion either to the right (as drawn) or to
the left
($l$ is odd in this picture).}
\label{figzigzag}
\end{figure}
%


The same reasoning gives the inequality
%
\begin{eqnarray}
\label{ineqkk1}
\hspace*{20pt}&&\P (LT_l,EU \mid Y_1^{m-1} \vee\bar Y_m^n )\nonumber
\\[-8pt]
\\[-8pt]
\hspace*{20pt}&& \qquad  \leq \biggl(\frac{p}{1-p} \biggr)^{u(l+1)+5m}  \biggl(\frac
12 \biggr)^{t(l+1)-t(l)}
\P (LT_{l+1},EU \mid Y_1^{m-1} \vee\bar Y_m^n ).
\nonumber
\end{eqnarray}
Indeed, any walk in $LT_l \cap EU$ covers $t(l+1)-t(l)$ more sites than
any walk in
$LT_{l+1} \cap EU$, while it is not hard to see that it makes at most
$u(l+1)+5m$ more
steps to the right. Since $t(l+1)-t(l) \sim4l$ and $u(l+1)+5m \sim2l$
as $l\to\infty$,
and $p/(1-p)<4$ (because $p<\frac45$), we find that
\[
\frac{
\P(LT_l,EU \mid Y_1^{m-1} \vee\bar Y_m^n)
}{
\P(LT_{l+1},EU \mid Y_1^{m-1} \vee\bar Y_m^n)
}
\]
decreases exponentially in $l$ for $l$ large. Hence the largest value
$l=L$ dominates.
Similar estimates allow us to neglect probabilities containing the
event~$NT$.

Combining (\ref{1split})--(\ref{ineqkk1}), we obtain that, for fixed $m$,
\[
\lim_{L\to\infty}\P (LT_L,EU \mid Y_1^{m-1} \vee\bar Y_m^n
) = 1,
\]
which immediately yields that, for fixed $m$,
%
\begin{equation}
\label{limids1}
 \P (C_0\,{=}\,B \mid Y_1^{m-1}\,{\vee}\,\bar Y_m^n )\,{=}\,\P (C_0\,{=}\,B
 \mid LT_L,EU,Y_1^{m-1}\,{\vee}\,\bar Y_m^n )  [1\,{+}\,o(1)],\hspace*{-40pt}
\end{equation}
where the error $o(1)$ tends to zero as $L\to\infty$.

The key point of (\ref{limids1}) is that $LT_L,EU,Y_1^{m-1} \vee\bar
Y_m^n$ forces the
coloring around the origin to look like $\cdots \mathit{BBWWBBWWBB} \cdots$.
More specifically,
$LT_L,EU,\bar Y_m^n$ tells us the coloring on a large region relative
to $S_m$ and, after
this, $Y_1^{m-1}$ determines the walk from time $1$ to time $m$
(relative to $S_1$).
Since~$S_1$ is independent of $\{LT_L,EU,Y_1^{m-1} \vee\bar Y_m^n\}$,
we therefore have
%
\begin{equation}
\label{lim1}
\P (C_0=B \mid LT_L,EU,Y_1^{m-1} \vee\bar Y_m^n ) \in\{
p,1-p\}.
\end{equation}

Equations (\ref{limids1}) and (\ref{lim1}) tell us that for large $n$,
$\P (C_0=B \mid Y_1^{m-1} \vee\bar Y_m^n )$ will be very
close to
$p$ or $1-p$. The idea now will be to modify the extension far away
so that an ``opposite'' type of structure is forced upon us and
thereby reverse the $p$ and $1-p$ above.

2. We next move to the choice of $\widetilde Y_m^n$. We take
%
\begin{eqnarray}
\label{Yblocksalt}
\widetilde Y_m^n &:=& [\mathit{WWBB}]^m \mathit{WBB} [\mathit{WWBB}]^{2m} \mathit{WBB}
[\mathit{WWBB}]^{2m+1} \nonumber
\\[-8pt]
\\[-8pt]
&&{} \cdots \mathit{WBB} [\mathit{WWBB}]^{2m+L-1} [\mathit{WWBB}]^{2m+L}.
\nonumber
\end{eqnarray}
The difference with $\bar Y_m^n$ in (\ref{Yblocks}) is that we removed
the last pivot $W_L$ and the 2 B's following it
(so that $n \to n-3$). The same computations as before give
%
\begin{eqnarray}
\label{limids2}
&&\P (C_0=B \mid Y_1^{m-1} \vee\widetilde Y_m^n )\nonumber
\\[-8pt]
\\[-8pt]
&& \qquad = \P (C_0=B \mid LT_{L-1},EU, Y_1^{m-1} \vee\widetilde Y_m^n
) [1+o(1)].
\nonumber
\end{eqnarray}
Now $LT_{L-1},EU, Y_1^{m-1} \vee\widetilde Y_m^n$ forces the walk to
do the exact opposite
up to time $t(L-1)$ to what $LT_L,EU, Y_1^{m-1} \vee\bar Y_m^n$ forced
it to do, because
there is one turn less and the walk still ends upwards. Therefore, by
symmetry, the walk
from time $1$ to time $m-1$ must also do the exact opposite, and so we
conclude that, for
$q\in\{p,1-p\}$,
%
\begin{eqnarray}
\label{symq}
&& \P (C_0=B \mid LT_L,EU,Y_1^{m-1} \vee\bar Y_m^n ) = q\nonumber
\\[-8pt]
\\[-8pt]
   && \qquad \Longleftrightarrow\quad
\P (C_0=B \mid LT_{L-1},EU,Y_1^{m-1} \vee\widetilde Y_m^n ) = 1-q.
\nonumber
\end{eqnarray}

Combining (\ref{limids1}) and (\ref{limids2})--(\ref{symq}), we obtain
the claim in
(\ref{limbartil}).
\end{pf}


\section{\texorpdfstring{$\B\notin\{\varnothing,\Omega\}$ for $p=\frac12$ and $\varepsilon=0$}
{B notin \{nothing,Omega\} for p = 1/2 and epsilon = 0}}
\label{S5}

In this section, we prove Theorem~\ref{phtr}(iii). We will prove that if $p=\frac12$ and $\eps=0$, then
%
\begin{eqnarray}
\label{badgood}
\hspace*{30pt}Y_1^\infty&=& B^\infty\ \mbox{is bad},\nonumber\\[-8pt]\\[-8pt]
\hspace*{30pt}Y_1^\infty&=& \mathit{BBWBB} [\mathit{WWBB}] \mathit{WBB} [\mathit{WWBB}]^2 \mathit{WBB} [\mathit{WWBB}]^3 \cdots\ \mbox{is good}.
\nonumber
\end{eqnarray}
(In the second line, $\mathit{BB}$ is put at the beginning to ensure that the
first $W$ may be a pivot.)


\subsection{\texorpdfstring{Proof of the first claim in (\protect\ref{badgood})}
{Proof of the first claim in (4.1)}}
\label{S5.1}
In this subsection, $Y_1^n$ and $Y_0^{n-1}$ denote random sequences,
and we switch back to specific sequences only in the last display.

Write
\begin{eqnarray*}
 \P (C_0 = W \mid Y_1^n = B^n )
&=& \P (C_0 = W \mid S_1=1, Y_1^n = B^n )\\
& =& \P (C_{-1}=W \mid Y_0^{n-1}=B^n ) = \frac{N(n)}{D(n)}
\end{eqnarray*}
with
%
\begin{eqnarray}
\label{NnDndef}
N(n) &:=& \P (C_{-1}=W,Y_0^{n-1}=B^n ) =
\sum_{i\in\N} \biggl(\frac12\biggr)^{i+2} p(n,i,1),\nonumber
\\[-8pt]
\\[-8pt]
D(n) &:=& \P (Y_0^{n-1}=B^n ) =
\sum_{i,j\in\N} \biggl(\frac12\biggr)^{i+j+1} p(n,i,j),
\nonumber
\end{eqnarray}
where $p(n,i,j) := \P(\tau_i \geq n,\tau_{-j}\geq n)$ is the
probability that simple
random walk (with $p=\frac12$ and $\eps=0$) starting from $0$ stays
between $-j+1$
and $i-1$ (inclusive) prior to time $n$. To see the second equality in
(\ref{NnDndef}),
let $E_{i,j}$ be the event that there is a $B$ at the origin, and the
first $W$ to
the right and to the left of the origin are located at $i$ and $-j$,
respectively. Then
\[
\P (Y_0^{n-1}=B^n )
= \sum_{i,j\in\N} \P(E_{i,j}) \P (Y_0^{n-1}=B^n\mid
E_{i,j} ),
\]
which is easily seen to be the claimed sum. The first equality in (\ref
{NnDndef}) is
handled similarly.

Trivially, $p(n,i,j) \geq p(n,i+j-1,1)$ for all $i,j\in\N$, and therefore
%
\begin{equation}
\label{Dnbd}
D(n) \geq\sum_{i\in\N} i\biggl(\frac12\biggr)^{i+2} p(n,i,1).
\end{equation}
Next, using Proposition 21.1 in \cite{S76}, we easily deduce that
\[
p(n,i,1) \sim \biggl[\cos\biggl(\frac{\pi}{i+1}\biggr) \biggr]^{n-1}
\cases{\displaystyle
C_i^\mathrm{even}, & \quad   as   $n \to\infty$  through   $n$   even,\cr\displaystyle
C_i^\mathrm{odd}, & \quad  as  $n \to\infty$  through   $n$   odd,
}
\]
where $\sim$ means that the ratio of the two sides tends to 1, and
\begin{eqnarray*}
C_i^\mathrm{even} &=& \frac{4}{i+1} \sin \biggl(\frac{\pi
}{i+1} \biggr)
\mathop{\mathop{\sum}_{ {0 \leq j < i} }}_{{j \mathrm{odd}} } \sin \biggl(\frac
{\pi(j+1)}{i+1} \biggr),\\
C_i^\mathrm{odd} &=& \frac{4}{i+1} \sin \biggl(\frac{\pi
}{i+1} \biggr)
\mathop{\mathop{\sum}_{ {0 \leq j < i} }}_{{j \mathrm{even}} } \sin \biggl(\frac
{\pi(j+1)}{i+1} \biggr).
\end{eqnarray*}
From this it follows that
%
\begin{equation}
\label{ratlim}
\lim_{n\to\infty} \frac{p(n,i+1,1)}{p(n,i,1)} = \infty, \qquad
i\in\N.
\end{equation}
Combining (\ref{NnDndef})--(\ref{ratlim}), we get $\lim_{n\to\infty}
N(n)/D(n)=0$, that is,
%
\begin{equation}
\label{limB1}
\lim_{n\to\infty} \P (C_0 = B \mid Y_1^n = B^n ) = 1.
\end{equation}

On the other hand, an extension of $Y_1^{m-1}=B^{m-1}$ with $\bar
Y_m^n$ as in Section~\ref{S4} gives
%
\begin{eqnarray}
\label{limB2}
\hspace*{20pt} P (C_0=B \mid Y_1^{m-1} \vee\bar Y_m^n )
&=& \P (C_0=B \mid LT_L, Y_1^{m-1} \vee\bar Y_m^n ) [1+o(1)]\nonumber
\\[-8pt]
\\[-8pt]
\hspace*{20pt} &=& \tfrac12 [1+o(1)]
\nonumber
\end{eqnarray}
[recall (\ref{ineqkk1})--(\ref{lim1})]. Combining
(\ref{limB1}) and (\ref{limB2}),
we get the first claim in (\ref{badgood}).\vspace*{6pt}


\subsection{\texorpdfstring{Proof of the second claim in (\protect\ref{badgood})}{Proof of the second claim in (4.1)}}
\label{S5.2}

Pick $L\in\N$ and $m-1=L(2L+ 5)+2$. Then
\[
Y_1^{m-1} = \mathit{BBWBB} [\mathit{WWBB}] \mathit{WBB} [\mathit{WWBB}]^2 \cdots \mathit{WBB} [\mathit{WWBB}]^L.
\]
As in Section \ref{S4}, a turn on a white pivot forces turns on all
previous white pivots. Therefore a walk compatible with $Y_1^{m-1}$
having at least one turn is characterized by the index
$k=0,1,\ldots,L-1$ of its last pivot $W_k$. The time of the $k$th
pivot is
$3+\sum_{j=0}^{k-1} [3+4(j+1)]$.

Conditioning on $Y_1^{m-1} \vee\bar Y_m^n$ still leaves us the freedom
to choose
$S_1\in\{-1,+1\}$ and $S_2 \in\{S_1-1,S_1+1\}$. Since $p=\frac12$,
it is easily
checked that, conditioned on $Y_1^{m-1} \vee\bar Y_m^n$ (and even on
the last pivot),
$S_1$ and $S_2-S_1$ are independent fair coin flips. There are $4$
compatible walks with
no turn and~$4L$ compatible walks with at least one turn. Since
$p=\frac12$, all
these walks have the same probability, but the walks with no turn have
a larger cost
for the coloring. Let $NT$ and $\mathit{AOT}:=[NT]^c$ denote the event that the
walk makes no turn,
respectively, at least one turn. We claim that
%
\begin{equation}
\label{Ptnt}
\P (NT \mid Y_1^{m-1} \vee\bar Y_m^n )
\leq\frac{1}{L+1} \P (\mathit{AOT} \mid Y_1^{m-1} \vee\bar Y_m^n ).
\end{equation}
To see how this comes about, recall Proposition \ref{decomp}, which
says that for an
arbitrary walk $s_1^{m-1}$ and an arbitrary extension $\bar Y_m^n$,
\begin{eqnarray*}
&&\P (S_1^{m-1}=s_1^{m-1}, Y_1^{m-1} \vee\bar Y_m^n )\\
&& \qquad  = \sum_{\bar s_1^{n-m+1}} \P(s_1^{m-1} \vee\bar
s_1^{n-m+1})
\biggl(\frac12\biggr)^{R(s_1^{m-1} \vee\bar s_1^{n-m+1})}\\
&&\hphantom{\sum_{\bar s_1^{n-m+1}}} \qquad  \quad {}\times
 1\{s_1^{m-1} \vee\bar s_1^{n-m+1} \sim Y_1^{m-1} \vee\bar Y_m^n\}.
\end{eqnarray*}
(The notation $s_1^{m-1} \vee\bar s_1^{n-m+1}$ denotes the walk
obtained by appending
the second walk to the end of the first walk.) Note that any compatible
walk up to time
$m-1$ ends either at the right end of the range or at the left end of
the range. Let
$s_1^{m-1}[0]$ and $s_1^{m-1}[1]$ denote compatible walks with no turn,
respectively, at
least one turn, either both ending at the right end of the range or
both ending at the
left end of the range. Then $R(s_1^{m-1}[0] \vee\bar s_1^{n-m+1}) \geq
R(s_1^{m-1}[1]
\vee\bar s_1^{n-m+1})$. Moreover, for any $\bar s_1^{n-m+1}$ and $\bar
Y_m^n$, if
$s_1^{m-1}[0] \vee\bar s_1^{n-m+1} \sim Y_1^{m-1} \vee\bar Y_m^n$,
then also $s_1^{m-1}[1]
\vee\bar s_1^{n-m+1} \sim Y_1^{m-1} \vee\bar Y_m^n$. Hence,
\[
\P (s_1^{m-1}[0], Y_1^{m-1} \vee\bar Y_m^n )
\leq\P (s_1^{m-1}[1], Y_1^{m-1} \vee\bar Y_m^n ).
\]
Summing over $s_1^{m-1}[0]$ and $s_1^{m-1}[1]$, we obtain (\ref{Ptnt}).

Next, on the event $\mathit{AOT}$, $C_0=B$ is fully determined by $S_1$ and
$S_2$. Therefore,
by symmetry,
\[
\P (C_0=B \mid \mathit{AOT}, Y_1^{m-1} \vee\bar Y_m^n ) =\tfrac12.
\]
Hence, uniformly in $\bar Y_m^n$,
\begin{eqnarray*}
&&\P (C_0=B \mid Y_1^{m-1} \vee\bar Y_m^n )\\
&& \qquad  = \P (C_0=B, \mathit{AOT} \mid Y_1^{m-1} \vee\bar Y_m^n )
+ \P (C_0=B, NT \mid Y_1^{m-1} \vee\bar Y_m^n )\\
&& \qquad  = \frac12 + O\biggl(\frac{1}{L}\biggr).
\end{eqnarray*}
Since $L\to\infty$ as $m\to\infty$, the second claim in
(\ref{badgood}) follows.\vspace*{6pt}


\section{\texorpdfstring{A possible approach to show that $\B=\Omega$ when $p \in[\frac12,\frac45)$ and $\eps\in(0,1)$}
{A possible approach to show that B = Omega when p in [1/2,4/5) and epsilon in (0,1)}}
\label{S3}

In this section, we explain a strategy for proving that $\B=\Omega$
when $p \in[\frac12,
\frac45)$ and $\eps\in(0,1)$. It seems that this case is much more
delicate than the
case $p \in[\frac12,\frac45)$ and $\eps= 0$ treated in Sections
\ref{S4}--\ref{S5}.
This strategy will be pursued in future work.


\subsection{Proposed strategy of the proof}
\label{S3.1}

For $M\in\N$, we use the notation $W^M$, $B^M$, $[\mathit{WB}]^M$ etc. to abbreviate
\[
\underbrace{\mathit{WW} \cdots W}_{M \mathrm{\ times\ } W} , \qquad\underbrace{\mathit{BB}
\cdots B}_{M \mathrm{\ times\ }B} ,
\qquad\underbrace{\mathit{WBWB}\cdots \mathit{WB}}_{M \mathrm{\ times\ }\mathit{WB}} ,\qquad\mbox{ etc.}
\]
Fix any configuration $Y_1^{\infty}$. To try to prove that $Y_1^\infty
$ is bad, we do the following:

(1) For $m,k,K\in\N$ with $k\geq2$, we consider the two color
records from time $m$ to time
$m+kK$ defined by
\[
\bar Y_m^{m+kK}(B) := [\mathit{WB}^{k-1}]^KW, \qquad\bar Y_m^{m+kK}(W) := [\mathit{BW}^{k-1}]^KB.
\]
(2) We \textit{expect} that, for any $p\in[\frac12,\frac45)$ and
$\eps\in(0,1)$,
%
\begin{eqnarray}
\label{star}
&&\inf_{m\in\N} \inf_{Y_1^{m-1}} \liminf_{k\to\infty} \liminf
_{K\to\infty}   \,\bigl|\P \bigl(C_0=B \mid Y_1^{m-1} \vee\bar Y_m^{m+kK}(B)
\bigr)\nonumber\\
&&\hphantom{\inf_{m\in\N} \inf_{Y_1^{m-1}} \liminf_{k\to\infty} \liminf
_{K\to\infty}   \bigl|}{}- \P\bigl (C_0=B \mid Y_1^{m-1} \vee\bar Y_m^{m+kK}(W) \bigr) \bigr|\\
 &&\qquad   \geq(1-\eps)(1-p).\nonumber
\end{eqnarray}
 In view of
Definition \ref{badconf}, this would imply
that $Y_1^\infty$ is a bad configuration, as desired.

The idea behind the above strategy is that $\bar Y_m^{m+kK}(B)$ forces
the walk to hit many
white sites at sparse times from time $m$ onwards. In order to achieve
this, the walk can
either move out to infinity, in which case the coloring must contain
many long black intervals,
or the walk can hang around the origin, in which case the coloring must
contain a single
white site close to the origin with two long black intervals on either
side. Since the
drift of the random walk is not too large, the best option is to hang
around the origin.
The single white site, at or next to the origin, is enough for the walk
to generate any~(!) color record $Y_1^{m-1}$ prior to time $m$, because the pausing
probability is strictly
positive. As a result, the conditional probability to see a black
origin given $Y_1^{m-1}
\vee\bar Y_m^{m+kK}(B)$ is closer to $1$ than given $Y_1^{m-1}\vee
\bar Y_m^{m+kK}(W)$. With
the latter conditioning, the role of $B$ and $W$ is reversed, and the
effect of the
conditioning is to have the origin lie in a region containing a single
black site separating
two long white intervals, so that the conditional probability to see a
black origin is closer to $0$.\vspace*{12pt}


\subsection{A few more details}
\label{S3.2}

The task is to control the conditional probability $\P(C_0=B \mid
Y_1^{m-1}\vee\bar
Y_m^{m+kK}(B))$. For that purpose, mark the positions of the walk at
the times $m+ki$,
$i=0,\ldots,K$, that correspond to the isolated~$W$'s in $\bar
Y_m^{m+kK}(B)$. By the
definition of $\bar Y_m^{m+kK}(B)$, two subsequent~$W$'s either
correspond to the same
white site or to two white sites that are separated by a single
interval of black sites of length at least~$1$.\looseness=1

\begin{figure}

\includegraphics{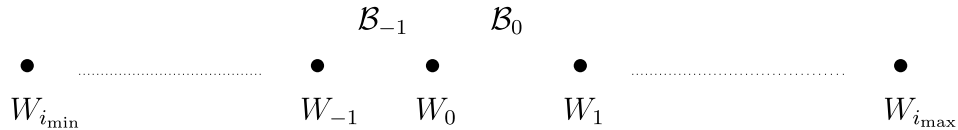}

\caption{White sites separated by black intervals. $W_0$
is the white site
seen at time $m$ in $Y_1^{m-1}\vee\bar Y_1^{m+kK}(B)$.}
\label{figcol}
\end{figure}

On the event $Y_1^{m-1}\vee\bar Y_m^{m+kK}(B)$, let $W_0$ be the white
site visited at
time~$m$. Relative to this site, all the white sites in $C$ can be
labeled $(W_i)_{i\in\Z}$,
with~$W_{-1}$ the first white site on the left of $W_0$, $W_1$ the
first white site on
the right of $W_0$, etc. (see Figure \ref{figcol}). Let $\cB_i$
denote the black interval
between~$W_i$ and $W_{i+1}$. $i_{\operatorname{min}}$ and $i_{\operatorname{max}}$ are the indices of
the left-most and
right-most white sites visited by the walk between times $m$ and $m+kK$.


The above representation allows to obtain an explicit (although complex)
formula for the conditional probability
$\P(\cdot  \mid   Y_1^{m-1} \vee Y_m^{m+kK} (B))$
involving classical simple random walk quantities.

Let $\cE_i$ denote the event that $\cB_i$ is visited between times
$m$ and
$m+kK$. Then the key fact that needs to be proved is the following:
%
\begin{equation}
\label{keyfact}
\hspace*{30pt}\inf_{Y_1^{m-1}} \liminf_{k\to\infty} \liminf_{K\to\infty}
\P \bigl(\cE_{-1} \cap\cE_0 \mid Y_1^{m-1} \vee\bar
Y_m^{m+kK}(B) \bigr) = 1
\qquad\forall m\in\N.
\end{equation}
From (\ref{keyfact}), we are able to prove the desired result (\ref{star}),
but the argument needed to prove (\ref{keyfact}) is long and we
are still working on trying to complete it.


\section*{Acknowledgments}
S\'{e}bastien Blach\`{e}re acknowledges
travel support from the ESF-program ``Random Dynamics in Spatially Extended
Systems'' (2002--2007). Frank den Hollander and Jeffrey E. Steif are
grateful for hospitality at the
Mittag-Leffler Institute in Stockholm in the Spring of 2009, when a~part
of the research on this paper was carried out.
We finally thank the referee for a thorough reading and comments on the paper.


%

\printaddresses

\end{document}